# Transversal Surfaces of Timelike Ruled Surfaces in Minkowski 3-Space $IR_1^3$

**Mehmet Önder**
*Celal Bayar University, Faculty of Science and Arts, Department of Mathematics, Muradiye Campus, 45047, Muradiye, Manisa, Turkey. E-mail: mehmet.onder@bayar.edu.tr*

**Abstract**
In this study we give definitions and characterizations of transversal surfaces of timelike ruled surfaces. We study some special cases such as the striction curve is a geodesic, an asymptotic line or a line of curvature. Moreover, we obtain developable conditions for transversal surfaces of a timelike ruled surface. In the paper, we consider timelike ruled surfaces with timelike rulings and timelike striction curve. The obtained results can be easily transferred to timelike ruled surfaces with spacelike rulings or to spacelike ruled surface.

## 1. Introduction

The notion of transversal surface was given by Sachs in 3-dimensional Euclidean space $E^3$ [6]. Sachs studied $\gamma$-transversal surface and new transversal surfaces called $\alpha$- and $\beta$-transversal surfaces of ruled surfaces in the same space. Taouktsoglou considered the notion of transversal surfaces for the ruled surfaces of the most general type in simply isotropic space $I_3^1$ [8]. Moreover, Sipus and Divjak described the transversal surfaces of ruled surfaces in the pseudo-Galilean space $G_3^1$ [7].

Let $IR_1^3$ be the Minkowski 3-space with standard Lorentzian flat metric given by
$$\langle , \rangle = -dx_1^2 + dx_2^2 + dx_3^2$$
where $(x_1, x_2, x_3)$ is a standard rectangular coordinate system of $IR_1^3$. An arbitrary vector $\vec{v} = (v_1, v_2, v_3)$ of $IR_1^3$ is said to be spacelike if $\langle \vec{v}, \vec{v} \rangle > 0$ or $\vec{v} = 0$, timelike if $\langle \vec{v}, \vec{v} \rangle < 0$ and null (lightlike) if $\langle \vec{v}, \vec{v} \rangle = 0$ and $\vec{v} \neq 0$. Similarly, an arbitrary curve $\vec{\alpha} = \vec{\alpha}(s)$ is said to be spacelike, timelike or null (lightlike), if all of its velocity vectors $\vec{\alpha}'(s)$ are spacelike, timelike or null (lightlike), respectively [5]. We say that a timelike vector is future pointing or past pointing if the first compound of vector is positive or negative, respectively. The norm of the vector $\vec{v} = (v_1, v_2, v_3)$ is given by $\|\vec{v}\| = \sqrt{|\langle \vec{v}, \vec{v} \rangle|}$.

For any vectors $\vec{x} = (x_1, x_2, x_3)$ and $\vec{y} = (y_1, y_2, y_3)$ in $IR_1^3$, Lorentzian vector product of $\vec{x}$ and $\vec{y}$ is defined by
$$\vec{x} \times \vec{y} = \begin{vmatrix} e_1 & -e_2 & -e_3 \\ x_1 & x_2 & x_3 \\ y_1 & y_2 & y_3 \end{vmatrix} = (x_2 y_3 - x_3 y_2, x_1 y_3 - x_3 y_1, x_2 y_1 - x_1 y_2).$$
(See [11]).

**Definition 1.1.(See [2,4])** *i) Hyperbolic angle:* Let $\vec{x}$ and $\vec{y}$ be future pointing (or past pointing) timelike vectors in $IR_1^3$. Then there is a unique real number $\theta \geq 0$ such that $<\vec{x},\vec{y}>=-\|\vec{x}\|\|\vec{y}\|\cosh\theta$. This number is called the *hyperbolic angle* between the vectors $\vec{x}$ and $\vec{y}$.

**ii) Central angle:** Let $\vec{x}$ and $\vec{y}$ be spacelike vectors in $IR_1^3$ that span a timelike vector subspace. Then there is a unique real number $\theta \geq 0$ such that $<\vec{x},\vec{y}>=\|\vec{x}\|\|\vec{y}\|\cosh\theta$. This number is called the *central angle* between the vectors $\vec{x}$ and $\vec{y}$.

**iii) Spacelike angle:** Let $\vec{x}$ and $\vec{y}$ be spacelike vectors in $IR_1^3$ that span a spacelike vector subspace. Then there is a unique real number $\theta \geq 0$ such that $<\vec{x},\vec{y}>=\|\vec{x}\|\|\vec{y}\|\cos\theta$. This number is called the *spacelike angle* between the vectors $\vec{x}$ and $\vec{y}$.

**iv) Lorentzian timelike angle:** Let $\vec{x}$ be a spacelike vector and $\vec{y}$ be a timelike vector in $IR_1^3$. Then there is a unique real number $\theta \geq 0$ such that $<\vec{x},\vec{y}>=\|\vec{x}\|\|\vec{y}\|\sinh\theta$. This number is called the *Lorentzian timelike angle* between the vectors $\vec{x}$ and $\vec{y}$.

**Definition 1.2.** A surface in the Minkowski 3-space $IR_1^3$ is called a timelike surface if the induced metric on the surface is a Lorentz metric and is called a spacelike surface if the induced metric on the surface is a positive definite Riemannian metric, i.e., the normal vector on the spacelike (timelike) surface is a timelike (spacelike) vector [3].

## 2. Timelike Ruled Surfaces in Minkowski 3-space

Let $I$ be an open interval in the real line $IR$. Let $\vec{f}=\vec{f}(u)$ be a curve in $IR_1^3$ defined on $I$ and $\vec{q}=\vec{q}(u)$ be a unit direction vector of an oriented line in $IR_1^3$. Then we have following parametrization for a ruled surface $N$

$$\vec{r}(u,v) = \vec{f}(u) + v\vec{q}(u). \qquad (1)$$

Assume that the surface normal is spacelike. Then by Definition 1.2, $N$ is a timelike ruled surface. The curve $\vec{f}=\vec{f}(u)$ is called base curve or generating curve of the surface and various positions of the generating lines $\vec{q}=\vec{q}(u)$ are called rulings. In particular, if the direction of $\vec{q}$ is constant, then the ruled surface is said to be cylindrical, and non-cylindrical otherwise.

The distribution parameter (or drall) of the timelike ruled surface in (1) is given by

$$d = \frac{\left|\dot{\vec{f}},\vec{q},\dot{\vec{q}}\right|}{\left\langle\dot{\vec{q}},\dot{\vec{q}}\right\rangle}, \qquad (2)$$

where $\dot{\vec{f}} = \dfrac{d\vec{f}}{du}$, $\dot{\vec{q}} = \dfrac{d\vec{q}}{du}$ (see [1,9]). If $\left|\dot{\vec{f}},\vec{q},\dot{\vec{q}}\right| = 0$, then normal vectors are collinear at all points of same ruling and at nonsingular points of the surface $N$, the tangent planes are identical. We then say that tangent plane contacts the surface along a ruling. Such a ruling is called a torsal ruling. If $\left|\dot{\vec{f}},\vec{q},\dot{\vec{q}}\right| \neq 0$, then the tangent planes of the surface $N$ are distinct at all points of same ruling which is called nontorsal [10].

**Definition 2.1.** A timelike ruled surface whose all rulings are torsal is called a developable timelike ruled surface. The remaining timelike ruled surfaces are called skew timelike ruled surfaces. Then, from (2) it is clear that a timelike ruled surface is developable if and only if at all its points the distribution parameter $d = 0$ [10].

For the unit normal vector $\vec{m}$ of a timelike ruled surface we have

$$\vec{m} = \frac{\vec{r}_u \times \vec{r}_v}{\|\vec{r}_u \times \vec{r}_v\|} = \frac{(\dot{\vec{f}} + v\dot{\vec{q}}) \times \vec{q}}{\sqrt{\left\langle \dot{\vec{f}}, \vec{q} \right\rangle^2 - \left\langle \vec{q}, \vec{q} \right\rangle \left\langle \dot{\vec{f}} + v\dot{\vec{q}},\ \dot{\vec{f}} + v\dot{\vec{q}} \right\rangle}}. \quad (3)$$

The unit normal of the surface along a ruling $u = u_1$ approaches a limiting direction as $v$ infinitely decreases. This direction is called the asymptotic normal direction and from (3) defined by

$$\vec{a} = \lim_{v \to \infty} \vec{m}(u_1, v) = \frac{\dot{\vec{q}} \times \vec{q}}{\|\dot{\vec{q}}\|}.$$

The point at which the unit normal of $N$ is perpendicular to $\vec{a}$ is called the striction point (or central point) $C$ and the set of striction points of all rulings is called striction curve of the surface. The parametrization of the striction curve $\vec{c} = \vec{c}(u)$ on a timelike ruled surface is given by

$$\vec{c}(u) = \vec{f}(u) + v_0 \vec{q}(u) = \vec{f} - \frac{\left\langle \dot{\vec{q}}, \dot{\vec{f}} \right\rangle}{\left\langle \dot{\vec{q}}, \dot{\vec{q}} \right\rangle} \vec{q}, \quad (4)$$

where $v_0 = -\dfrac{\left\langle \dot{\vec{q}}, \dot{\vec{f}} \right\rangle}{\left\langle \dot{\vec{q}}, \dot{\vec{q}} \right\rangle}$ is called strictional distance [10].

The vector $\vec{h}$ defined by $\vec{h} = \pm \vec{a} \times \vec{q}$ is called central normal which is the surface normal along the striction curve. Then the orthonormal system $\{C; \vec{q}, \vec{h}, \vec{a}\}$ is called Frenet frame of the ruled surfaces $N$ where $C$ is the central point of ruling of timelike ruled surface $N$ and $\vec{q}, \vec{h} = \pm \vec{a} \times \vec{q}, \vec{a}$ are unit vectors of ruling, central normal and central tangent, respectively.

Let now consider ruled surface $N$ with non-null Frenet vectors and their non-null derivatives. According to the Lorentzian character of ruling, we can give the following classifications of the timelike ruled surface $N$ as follows;

**i)** If $\vec{q}$ is timelike, then timelike ruled surface $N$ is said to be of type $N_-$.

**ii)** If $\vec{q}$ is spacelike, then timelike ruled surface $N$ is said to be of type $N_+$.

In these classifications we use subscript "+" and "-" to show the Lorentzian casual character of ruling. By using these classifications, parametrization of timelike ruled surface $N$ can be given as follows,

$$\vec{r}(u,v) = \vec{f}(u) + v \vec{q}(u),$$

where $\langle \vec{q}, \vec{q} \rangle = \varepsilon (= \pm 1)$, $\langle \vec{h}, \vec{h} \rangle = 1$, $\langle \vec{a}, \vec{a} \rangle = -\varepsilon$.

For the derivatives of vectors of Frenet frame $\{C; \vec{q}, \vec{h}, \vec{a}\}$ of timelike ruled surface $N$ with respect to the arc length $s$ of striction curve we have

$$\begin{bmatrix} d\vec{q}/ds \\ d\vec{h}/ds \\ d\vec{a}/ds \end{bmatrix} = \begin{bmatrix} 0 & k_1 & 0 \\ -\varepsilon k_1 & 0 & k_2 \\ 0 & \varepsilon k_2 & 0 \end{bmatrix} \begin{bmatrix} \vec{q} \\ \vec{h} \\ \vec{a} \end{bmatrix} \quad (5)$$

where $k_1 = \dfrac{ds_1}{ds}$, $k_2 = \dfrac{ds_3}{ds}$ and $s_1$, $s_3$ are the arc lengths of the spherical curves circumscribed by the bound vectors $\vec{q}$ and $\vec{a}$, respectively. Timelike ruled surfaces satisfying $k_1 \neq 0$, $k_2 = 0$ are called timelike conoids (For details [10]).

In this study we introduce the definitions and characterizations of transversal surfaces of timelike ruled surfaces. For a reference ruled surface, we consider a timelike ruled surface of the type $N_-$ with timelike striction curve. So, we first give some special cases for the striction curve of a timelike ruled surface of the type $N_-$. Of course the obtained results of the following sections can be easily transferred to other cases such as the surface is of the type $N_+$ or is a spacelike ruled surface.

## 3. Some Special Cases for the Striction Curve of a Timelike Ruled Surface

Let assume that timelike ruled surface $N$ be of the type $N_-$ and let the generating curve of the surface be its striction curve and the Frenet frame of the surface be $\{C; \vec{q}, \vec{h}, \vec{a}\}$. Moreover, assume that the striction curve is timelike. Then for the parametrization of the surface $N$ we write

$$\vec{r}(s,v) = \vec{c}(s) + v\vec{q}(s), \tag{6}$$

where $\langle \vec{q}, \vec{q} \rangle = -1$, $\langle \vec{h}, \vec{h} \rangle = 1$, $\langle \vec{a}, \vec{a} \rangle = 1$ and $s$ is the arc length parameter of striction curve $\vec{c}(s)$. For the unit tangent vector of the striction curve we have

$$\vec{c}'(s) = \dfrac{d\vec{c}}{ds} = \vec{t} = \cosh\theta\, \vec{q} + \sinh\theta\, \vec{a}, \tag{7}$$

where $\theta$ is the hyperbolic angle between timelike vectors $\vec{t}$ and $\vec{q}$. Then from (2) and (7) the distribution parameter of the surface is obtained as

$$d = -\dfrac{\sinh\theta}{k_1}. \tag{8}$$

Let now investigate some special cases for the striction curve $\vec{c}(s)$. First assume that striction curve $\vec{c}(s)$ be an asymptotic line on $N$. Then central normal $\vec{h}$ and direction vector $\vec{c}''$ of principal normal vector satisfy $\langle \vec{h}, \vec{c}'' \rangle = 0$. After a simple computation we have $\tanh\theta = k_1/k_2$. Then we have the following theorem.

***Theorem 3.1.*** *Let timelike ruled surface $N$ be of the type $N_-$. Then timelike striction curve $\vec{c}(s)$ of the surface is an asymptotic line on $N$ if and only if $\tanh\theta = k_1/k_2$ holds.*

If the striction curve $\vec{c}(s)$ is a geodesic on $N$, then central normal $\vec{h}$ and direction vector $\vec{c}''$ of principal normal vector are linearly depended i.e., we have $\vec{h} = \lambda \vec{c}''$ where $\lambda = \lambda(s)$ is a scalar function. The last equality gives us

$$\vec{h} = \lambda\big((\theta'\sinh\theta)\vec{q} + (k_1\cosh\theta - k_2\sinh\theta)\vec{h} + (\theta'\cosh\theta)\vec{a}\big), \tag{9}$$

and from (9) we have that $\theta$ is constant and we give the following theorem.

**Theorem 3.2.** *Let timelike ruled surface $N$ be of the type $N_-$. Then timelike striction curve $\vec{c}(s)$ of the surface is a geodesic on $N$ if and only if $\theta$ is constant.*

Finally, assume that the striction curve $\vec{c}(s)$ is a line of curvature on $N$. Then the derivative of central normal $\vec{h}$ and tangent vector $\vec{t}$ of striction curve $\vec{c}(s)$ are linearly depended i.e., we have $\vec{h}' = \lambda \vec{t}$ where $\lambda = \lambda(s)$ is a scalar function. Then from (5) and (7) we have $\tanh\theta = k_2/k_1$ which gives us following theorem.

**Theorem 3.3.** *Let timelike ruled surface $N$ be of the type $N_-$. Then timelike striction curve $\vec{c}(s)$ of the surface is a line of curvature on $N$ if and only if $\tanh\theta = k_2/k_1$ holds.*

Now we can introduce transversal surfaces of a timelike ruled surface of the type $N_-$. When we talk about timelike ruled surface $N$ and striction curve $\vec{c}(s)$, we mean that surface is of the type $N_-$ and striction curve $\vec{c}(s)$ is timelike and for short we don't write the Lorentzian characters of the surface and curve.

### 4. $\alpha$-Transversal Surfaces of Timelike Ruled Surfaces

In this section we give the definition and characterizations of $\alpha$-transversal surfaces of a timelike ruled surface. First, we give the following definition.

**Definition 4.1.** Let $N$ be a timelike ruled surface. An $\alpha$-transversal surface $N^\alpha$ of $N$ is a ruled surface in $IR_1^3$ whose rulings are straight lines through a striction point $\vec{c}(s)$ determined by ruling $\vec{q}_\alpha = \mu(\alpha)\vec{q} + \eta(\alpha)\vec{h}$ where

$$\mu(\alpha) = \begin{cases} \cosh\alpha, & \vec{q}_\alpha \text{ is timelike} \\ \sinh\alpha, & \vec{q}_\alpha \text{ is spacelike} \end{cases} \quad \eta(\alpha) = \begin{cases} \sinh\alpha, & \vec{q}_\alpha \text{ is timelike} \\ \cosh\alpha, & \vec{q}_\alpha \text{ is spacelike} \end{cases} \tag{10}$$

and to obtained the non-trivial cases (ruling is not $\vec{q}$ or $\vec{h}$) we assume $\mu(\alpha) \neq 0$, $\eta(\alpha) \neq 0$.

From Definition 4.1. the parametrization of $\alpha$-transversal surface $N^\alpha$ is
$$\vec{r}_\alpha(s,v) = \vec{c}(s) + v\vec{q}_\alpha(s), \tag{11}$$
and from (10) we have
$$\langle \vec{q}_\alpha, \vec{q}_\alpha \rangle = \eta^2 - \mu^2 = \ell = \pm 1. \tag{12}$$
The strictional distance $v_\alpha$ of the $\alpha$-transversal surface $N^\alpha$ is obtained as
$$v_\alpha = -\frac{\langle \vec{c}', \vec{q}_\alpha' \rangle}{\langle \vec{q}_\alpha', \vec{q}_\alpha' \rangle} = \frac{\eta(\cosh\theta(\alpha' + k_1) - \sinh\theta k_2)}{\eta^2 k_2^2 - \ell(\alpha' + k_1)^2}, \tag{13}$$
where $\alpha' = d\alpha/ds$. Then we obtain the following theorem:

**Theorem 4.1.** *The striction curve $\vec{c}_\alpha$ on every $N^\alpha$ coincides with the striction curve $\vec{c}(s)$ if and only if*
$$\tanh\theta = \frac{\alpha' + k_1}{k_2}, \tag{14}$$
*holds.*

From Theorem 3.1 we know that the striction curve $\vec{c}(s)$ is an asymptotic line on $N$ if and only if $\tanh\theta = k_1/k_2$ holds. In this special case, Theorem 4.1 gives us $\alpha' = 0$. Then we can give the following theorem:

**Theorem 4.2.** *Let the striction curve $\vec{c}_\alpha$ on every $N^\alpha$ coincides with the striction curve $\vec{c}(s)$. Then $\vec{c}(s)$ is an asymptotic line on $N$ if and only if $\alpha$ is constant.*

We know that striction curve $\vec{c}(s)$ of the surface is a geodesic on $N$ if and only if $\theta$ is constant. Then from (14) we have the following theorem:

**Theorem 4.3.** *Let the striction curve $\vec{c}_\alpha$ on every $N^\alpha$ coincides with the striction curve $\vec{c}(s)$. Then $\vec{c}(s)$ is a geodesic on $N$ if and only if there exists a constant $x$ such that $\alpha' = xk_2 - k_1$.*

From Theorem 3.3, the striction curve $\vec{c}(s)$ is a line of curvature on $N$ if and only if $\tanh\theta = k_2/k_1$ holds. Then from (14) we have the following theorem:

**Theorem 4.4.** *Let the striction curve $\vec{c}_\alpha$ on every $N^\alpha$ coincides with the striction curve $\vec{c}(s)$. Then $\vec{c}(s)$ is a line of curvature on $N$ if and only if $\alpha' = \dfrac{k_2^2 - k_1^2}{k_1}$ holds.*

Let now consider the developable $\alpha$-transversal surfaces. By a simple calculation from (2) and (11) the distribution parameter $d_\alpha$ of $N^\alpha$ is obtained as

$$d_\alpha = \frac{\ell(\alpha' + k_1)\sinh\theta - \eta^2 k_2 \cosh\theta}{\eta^2 k_2^2 - \ell(\alpha' + k_1)^2}. \tag{15}$$

Then we have the following theorem:

**Theorem 4.5.** *$\alpha$-transversal surface $N^\alpha$ is developable if and only if $\tanh\theta = \dfrac{\ell(\alpha' + k_1)}{\mu^2 k_2}$ holds.*

Moreover, from (8) and (15) we have

$$d_\alpha = -\frac{\ell d k_1(\alpha' + k_1) + \eta^2 k_2 \cosh\theta}{\eta^2 k_2^2 - \ell(\alpha' + k_1)^2}. \tag{16}$$

where $d$ is distribution parameter of $N$. Since we consider non-trivial cases i.e., $\mu \neq 0$, if $d = 0$ then from (16) we have $k_2 = 0$ which gives us following corollary:

**Corollary 4.6.** *Let timelike ruled surface $N$ be developable. Then $N^\alpha$ is developable if and only if $N$ is a timelike conoid.*

## 5. $\beta$-Transversal Surfaces of Timelike Ruled Surfaces

In this section we give the definition and characterizations of $\beta$-transversal surfaces of a timelike ruled surface. First, we give the following definition.

**Definition 5.1.** Let $N$ be a timelike ruled surface in $IR_1^3$. The $\beta$-transversal surface $N^\beta$ of $N$ is a ruled surface in $IR_1^3$ whose rulings are straight lines through a striction point $\vec{c}(s)$ determined by ruling $\vec{q}_\beta = \cos\beta \vec{h} + \sin\beta \vec{a}$ where $\beta$ is spacelike angle between $\vec{q}_\beta$ and $\vec{h}$, and to obtained the non-trivial cases (ruling is not $\vec{h}$ or $\vec{a}$) we assume $\beta \neq n\pi$, $(2n+1)\pi/2$ where $n \in \mathbb{Z}$.

From this definition the parametrization of $\beta$-transversal surface $N^\beta$ is
$$\vec{r}_\beta(s,v) = \vec{c}(s) + v\vec{q}_\beta(s). \tag{17}$$
By a simple calculation the strictional distance $v_\beta$ of the $\beta$-transversal surface $N^\beta$ is obtained as
$$v_\beta = -\frac{\langle \vec{c}', \vec{q}'_\beta \rangle}{\langle \vec{q}'_\beta, \vec{q}'_\beta \rangle} = \frac{\cos\beta((\beta' + k_2)\sinh\theta - k_1 \cosh\theta)}{(\beta' + k_2)^2 - k_1^2 \cos^2\beta}, \tag{18}$$
where $\beta' = d\beta/ds$. Then we obtain the following theorem:

**Theorem 5.1.** *The striction curve $\vec{c}_\beta$ on every $N^\beta$ coincides with the striction curve $\vec{c}(s)$ if and only if*
$$\tanh\theta = \frac{k_1}{\beta' + k_2}, \tag{19}$$
*holds.*

From Theorem 3.1 we have that the striction curve $\vec{c}(s)$ is an asymptotic line on $N$ if and only if $\tanh\theta = k_1/k_2$ holds. In this special case, Theorem 5.1 gives us $\beta' = 0$. Then we can give the following theorem:

**Theorem 5.2.** *Let the striction curve $\vec{c}_\beta$ on every $N^\beta$ coincides with the striction curve $\vec{c}(s)$. Then $\vec{c}(s)$ is an asymptotic line on $N$ if and only if $\beta$ is constant.*

We know that striction curve $\vec{c}(s)$ of the surface is a geodesic on $N$ if and only if $\theta$ is constant (Theorem 3.2). Then from (19) we have the following theorem:

**Theorem 5.3.** *Let the striction curve $\vec{c}_\beta$ on every $N^\beta$ coincides with the striction curve $\vec{c}(s)$. Then $\vec{c}(s)$ is a geodesic on $N$ if and only if there exists a constant $y$ such that $k_1 = y(\beta' + k_2)$.*

From Theorem 3.3, the striction curve $\vec{c}(s)$ is a line of curvature on $N$ if and only if $\tanh\theta = k_2/k_1$ holds. Then from (19) we have the following theorem:

**Theorem 5.4.** *Let the striction curve $\vec{c}_\beta$ on every $N^\beta$ coincides with the striction curve $\vec{c}(s)$. Then $\vec{c}(s)$ is a curvature line on $N$ if and only if $\beta' = \frac{k_1^2 - k_2^2}{k_2}$ holds.*

Let now consider the developable $\beta$-transversal surfaces. By a simple calculation from (2) and (17) the distribution parameter $d_\beta$ of $N^\beta$ is obtained as

$$d_\beta = \frac{k_1 \cos^2 \beta \sinh \theta - (\beta' + k_2) \cosh \theta}{(\beta' + k_2)^2 - k_1^2 \cos^2 \beta}. \tag{20}$$

Then we have the following theorem:

**Theorem 5.5.** *$\beta$-transversal surface $N^\beta$ is developable if and only if $\tanh \theta = \dfrac{\beta' + k_2}{k_1 \cos^2 \beta}$ holds.*

Moreover, from (8) and (20) we have

$$d_\beta = -\frac{dk_1^2 \cos^2 \beta + (\beta' + k_2) \cosh \theta}{(\beta' + k_2)^2 - k_1^2 \cos^2 \beta}, \tag{21}$$

where $d$ is distribution parameter of reference surface $N$. Then (21) gives us following corollary:

**Corollary 5.6.** *Let timelike ruled surface $N$ be developable. Then $N^\beta$ is developable if and only if $\beta' = -k_2$.*

## 6. $\gamma$-Transversal Surfaces of Timelike Ruled Surfaces

In this section we give the definition and characterizations of $\gamma$-transversal surfaces of a timelike ruled surface. First, we give the following definition.

**Definition 6.1.** Let $N$ be a timelike ruled surface. The $\gamma$-transversal surface $N^\gamma$ of $N$ is a ruled surface in $IR_1^3$ whose rulings are straight lines through a striction point $\vec{c}(s)$ determined by ruling $\vec{q}_\gamma = \mu(\gamma)\vec{q} + \eta(\gamma)\vec{a}$ where

$$\mu(\gamma) = \begin{cases} \cosh \gamma, & \vec{q}_\gamma \text{ is timelike} \\ \sinh \gamma, & \vec{q}_\gamma \text{ is spacelike} \end{cases} \quad \eta(\gamma) = \begin{cases} \sinh \gamma, & \vec{q}_\gamma \text{ is timelike} \\ \cosh \gamma, & \vec{q}_\gamma \text{ is spacelike} \end{cases} \tag{22}$$

and to obtained the non-trivial cases (ruling is not $\vec{q}$ or $\vec{a}$) we assume $\mu(\gamma) \neq 0$, $\eta(\gamma) \neq 0$.

From Definition 6.1. the parametrization of $\gamma$-transversal surface $N^\gamma$ is

$$\vec{r}_\gamma(s, v) = \vec{c}(s) + v \vec{q}_\gamma(s), \tag{23}$$

and from (22) we have

$$\langle \vec{q}_\gamma, \vec{q}_\gamma \rangle = \eta^2 - \mu^2 = \ell = \pm 1. \tag{24}$$

The strictional distance $v_\gamma$ of the $\gamma$-transversal surface $N^\gamma$ is obtained as

$$v_\gamma = -\frac{\langle \vec{c}', \vec{q}_\gamma' \rangle}{\langle \vec{q}_\gamma', \vec{q}_\gamma' \rangle} = \frac{\gamma'(\mu \sinh \theta - \eta \cosh \theta)}{(\mu k_1 - \eta k_2)^2 - \ell(\gamma')^2}, \tag{25}$$

where $\gamma' = d\gamma/ds$. Then from (25) we obtain the following theorem:

**Theorem 6.1.** *The striction curve $\vec{c}_\gamma$ on every $N^\gamma$ coincides with the striction curve $\vec{c}(s)$ if and only if $\gamma$ is constant or*

$$\tanh\theta = \frac{\eta}{\mu}, \tag{26}$$

holds.

From Theorem 3.1 we know that the striction curve $\vec{c}(s)$ is an asymptotic line on $N$ if and only if $\tanh\theta = k_1/k_2$ holds. Then we can give the following theorem:

**Theorem 6.2.** *Let the striction curve $\vec{c}_\gamma$ on every $N^\gamma$ coincides with the striction curve $\vec{c}(s)$ and let $\gamma$ be non-constant. Then $\vec{c}(s)$ is an asymptotic line on $N$ if and only if $\dfrac{k_1}{k_2} = \dfrac{\eta}{\mu}$ holds.*

We know that striction curve $\vec{c}(s)$ of the surface is a geodesic on $N$ if and only if $\theta$ is constant. Then from (26) we have the following theorem:

**Theorem 6.3.** *Let the striction curve $\vec{c}_\gamma$ on every $N^\gamma$ coincides with the striction curve $\vec{c}(s)$ and let $\gamma$ be non-constant. Then $\vec{c}(s)$ is a geodesic on $N$ if and only if there exists a constant $z$ such that $\eta = z\mu$.*

From Theorem 3.3, the striction curve $\vec{c}(s)$ is a line of curvature on $N$ if and only if $\tanh\theta = k_2/k_1$ holds. Then from (26) we have the following theorem:

**Theorem 6.4.** *Let the striction curve $\vec{c}_\gamma$ on every $N^\gamma$ coincides with the striction curve $\vec{c}(s)$ and let $\gamma$ be non-constant. Then $\vec{c}(s)$ is a line of curvature on $N$ if and only if $\dfrac{k_2}{k_1} = \dfrac{\eta}{\mu}$ holds.*

Let now consider the developable $\gamma$-transversal surfaces. By a simple calculation from (2) and (23) the distribution parameter of $N^\gamma$ is obtained as

$$d_\gamma = \frac{(\eta\cosh\theta - \mu\sinh\theta)(\mu k_1 - \eta k_2)}{(\mu k_1 - \eta k_2)^2 - \ell(\gamma')^2}. \tag{27}$$

Then we have the following theorem:

**Theorem 6.5.** *$\gamma$-transversal surface $N^\gamma$ is developable if and only if $\tanh\theta = \dfrac{\eta}{\mu}$ or $\dfrac{k_1}{k_2} = \dfrac{\eta}{\mu}$ hold.*

Moreover, from (8) and (27) we have

$$d_\gamma = \frac{(\eta\cosh\theta + k_1 d\mu)(\mu k_1 - \eta k_2)}{(\mu k_1 - \eta k_2)^2 - \ell(\gamma')^2}. \tag{28}$$

where $d$ is distribution parameter of $N$. Then if $N$ is developable, i.e., $d = 0$, (28) gives us following corollary.

***Corollary 6.6.*** *Let timelike ruled surface* $N$ *be developable. Then* $N^\gamma$ *is developable if and only if* $\dfrac{k_1}{k_2} = \dfrac{\eta}{\mu}$ *holds.*

## 7. Conclusions

Transversal surfaces of a timelike ruled surface of the type $N_-$ are defined and characterizations of these surfaces are given. In the paper we consider the striction line as a timelike curve. Of course, one can obtain corresponding theorems for a timelike ruled surface with a spacelike striction line, for a surface of the type $N_+$ or for a spacelike ruled surface.

## References

[1] Abdel-All, N.H., Abdel-Baky, R.A., Hamdoon, F.M., *Ruled surfaces with timelike rulings,* App. Math. and Comp. 147, 241–253 (2004).

[2] Aydoğmuş, Ö., Kula, L., Yaylı, Y., *On point-line displacement in Minkowski 3-space,* Differ. Geom. Dyn. Syst. 10, 32-43 (2008).

[3] Beem, J.K., Ehrlich, P.E., Global Lorentzian Geometry, Marcel Dekker, New York, (1981).

[4] Birman, G., Nomizo, K., *Trigonometry in Lorentzian Geometry,* Ann. Math. Mont. 91, (9), 543-549 (1984).

[5] O'Neill, B., Semi-Riemannian Geometry with Applications to Relativity. Academic Press, London (1983).

[6] Sachs, H., *Uber Transversalflachen von Regelflachen,* Sitzungsber. Akad. Wiss. Wien, 186, 427–439 (1978).

[7] Sipus, Z.M., Divjak, B., *Transversal Surfaces of Ruled Surfaces in the Pseudo-Galilean Space,* Sitzungsber. Abt. II, 213, 23–32 (2004).

[8] Taouktsoglou, A., *Transversalflachen von Regelflachen im einfach isotropen Raum,* Sitzungsber. Akad. Wiss. Wien, 203, 137–148 (1994).

[9] Turgut, A., Hacısalihoğlu, H.H., *Timelike ruled surfaces in the Minkowski 3-space,* Far East J. Math. Sci. 5 (1), 83–90 (1997).

[10] Uğurlu, H.H., Önder, M., *Instantaneous Rotation vectors of Skew Timelike Ruled Surfaces in Minkowski 3-space,* VI. Geometry Symposium, Uludağ University, Bursa, Turkey. (2008).

[11] Walrave, J., *Curves and surfaces in Minkowski space,* PhD. thesis, K.U. Leuven, Fac. of Science, Leuven, (1995).